\documentclass[11pt,a4paper]{article}
\usepackage[margin=27mm]{geometry}
\usepackage{amsmath,amssymb,amsthm,mathtools,microtype}
\usepackage[colorlinks=true,linkcolor=blue,citecolor=blue,urlcolor=blue]{hyperref}
\allowdisplaybreaks[2]
\numberwithin{equation}{section}
\newtheorem{theorem}{Theorem}[section]
\newtheorem{lemma}[theorem]{Lemma}
\newtheorem{proposition}[theorem]{Proposition}
\theoremstyle{remark}

\newcommand{\R}{\mathbb R}
\newcommand{\E}{\mathbb E}
\newcommand{\Pp}{\mathbb P}
\newcommand{\one}{\mathbf1}
\newcommand{\Capacity}{\operatorname{Cap}}
\newcommand{\supp}{\operatorname{supp}}
\newcommand{\TV}{\operatorname{TV}}
\title{Polarity of the First Dirichlet Liouville Spectral Gap}
\author{Chunhao Cai\\
School of Mathematics (Zhuhai), Sun Yat-Sen University, China\\
\texttt{caichh9@mail.sysu.edu.cn}}
\date{}
\begin{document}
\maketitle

\begin{abstract}
We prove zero $1$-capacity of first-gap collisions for the minimal
Liouville spectral Green-response form and its canonical extension.
The minimal-form diffusion avoids collisions for every initial law
absolutely continuous with respect to the spectral law.
\end{abstract}



\section{Introduction}\label{sec:introduction}

The first Dirichlet Liouville eigenvalue is simple, so the first two
eigenvalues are strictly ordered for almost every realization of the
Gaussian free field.  This does not determine whether a diffusion of
the spectral coordinates can reach the diagonal where they coincide.
We study the capacity of this diagonal for the Dirichlet forms induced
by the Gaussian response of the Liouville spectrum.  The static
spectral and response results used below are established in
\cite{Cai,CaiMV}.

Liouville Brownian motion was constructed in \cite{GRV}, and its
spectral Weyl asymptotics were studied in \cite{BW}.  Capacity
arguments for noncollision in determinantal dynamics appear in
\cite{Osada}.  Here the spectral form is induced by Gaussian--GMC
response, with polarity proved through Green localization and
inverse-gap estimates.

Let $D\subset\R^2$ be smooth, bounded and connected.  Set
\[
 H=H_0^1(D),\qquad
 \langle u,v\rangle_H=\frac1{2\pi}\int_D\nabla u\cdot\nabla v\,dx.
\]
On a probability space $(\Omega,\mathcal F,\Pp)$, let $h$ be the
zero-boundary Gaussian free field on $D$ with covariance $G_D$,
normalized to have logarithmic singularity of coefficient one.  Fix
$0<\gamma<2$ and let $\mu=M_\gamma(h)$ be the corresponding
Wick-normalized Gaussian multiplicative chaos measure.  Write the
Dirichlet Liouville eigenpairs as
\begin{equation}\label{eq:eigenpairs}
 \langle\phi_i,u\rangle_H=\lambda_i\int_D\phi_i u\,d\mu,
 \qquad\int_D\phi_i\phi_j\,d\mu=\delta_{ij},
 \qquad0<\lambda_1<\lambda_2\le\cdots .
\end{equation}
For a fixed integer $N\ge2$, put
\begin{equation}\label{eq:spectral-data}
 Y_i=\log\lambda_i,\qquad
 Y=(Y_1,\ldots,Y_N),\qquad g=Y_2-Y_1,
 \qquad\rho=\operatorname{Law}(Y),\qquad S=\supp\rho.
\end{equation}
Static simplicity gives $\Pp(g>0)=1$ \cite{Cai}.  Hence
$\rho\{y\in S:y_1<y_2\}=1$ and
$S\subseteq\{y\in\R^N:y_1\le y_2\}$.  Define the collision set
\[
 \Sigma_1=\{y\in S:y_1=y_2\}.
\]

Let $\nabla_h$ denote the Gaussian gradient.  To define the spectral
response forms, let $-\Delta_D\psi_k=\kappa_k\psi_k$, where $(\psi_k)$
is orthonormal in $L^2(D,dx)$, and set
\begin{equation}\label{eq:H-basis}
 e_k=\sqrt{\frac{2\pi}{\kappa_k}}\,\psi_k,\qquad
 P_se_k=(1+\kappa_k)^{-s}e_k,\qquad s\ge0.
\end{equation}
On $C_c^\infty(\R^N)|_S$, define
\begin{equation}\label{eq:weighted-form}
 E_s(u,v)=\E\langle P_s\nabla_h(u(Y)),\nabla_h(v(Y))\rangle_H,
 \qquad s\ge0.
\end{equation}
For $s>0$, let $(E_s,F_s)$ be the minimal closure of this preform;
write $(E_0^{\min},F_0^{\min})$ for its minimal closure at $s=0$.
As $s\downarrow0$ the forms increase.  For any sequence
$s_j>0$ with $s_j\downarrow0$, put
\begin{equation}\label{eq:canonical-form-domain}
 \begin{split}
 F_*&=\left\{u\in\bigcap_jF_{s_j}:
                  \sup_jE_{s_j}(u,u)<\infty\right\},\\
 E_*(u,u)&=\sup_j E_{s_j}(u,u),\qquad u\in F_*.
 \end{split}
\end{equation}
By \cite[Proposition A.6 and Theorem A.18]{CaiMV}, this definition is
independent of $(s_j)$, the minimal form is regular, strongly local
and conservative, the canonical form is closed, strongly local and
conservative, and
\begin{equation}\label{eq:form-interface}
 F_0^{\min}\subseteq F_*,\qquad
 E_*(u,u)=E_0^{\min}(u,u)\quad(u\in F_0^{\min}).
\end{equation}

For a closed form $(E,F)$ on $L^2(S,\rho)$, its outer $1$-capacity is
\[
 \Capacity_{E,1}(A)
 =\inf_{\substack{O\supset A\\ O\text{ relatively open in }S}}
   \inf_{\substack{u\in F\\u\ge1\ \rho\text{-a.e. on }O}}
   \bigl(E(u,u)+\|u\|_{L^2(\rho)}^2\bigr).
\]
Let $X^{\min}$ be a conservative continuous Hunt process properly
associated with $(E_0^{\min},F_0^{\min})$ on $S$.  For any probability
measure $\pi\ll\rho$, denote its law with initial distribution $\pi$
by $\Pp_\pi$.

\begin{theorem}\label{thm:noncollision}
For every $0<\gamma<2$,
\begin{equation}\label{eq:main-capacity}
 \Capacity_{E_0^{\min},1}(\Sigma_1)
 =\Capacity_{E_*,1}(\Sigma_1)=0.
\end{equation}
Moreover, for every probability measure $\pi\ll\rho$,
\begin{equation}\label{eq:main-process}
 \Pp_\pi\{X^{\min}_1(t)<X^{\min}_2(t)
                \text{ for all }t\ge0\}=1.
\end{equation}
\end{theorem}

The proof first bounds the spectral gap by spatial localization and
then establishes inverse-gap moments along Gaussian rotations.  These
estimates yield a smoothed response-energy bound; a high-frequency
cutoff completes the capacity argument in
Section~\ref{sec:final-proof}.  The uniform GMC and rotation estimates
are proved in the appendices.

\section{Preliminaries}\label{sec:preliminaries}

In the basis \eqref{eq:H-basis}, the Gaussian free field has the
expansion $h=\sum_{k\ge1}\xi_ke_k$ in the distributional sense,
where $(\xi_k)$ are independent standard Gaussian variables.
The Gaussian gradient $\nabla_h$ is the closure of
$\sum_{k=1}^m(\partial_{\xi_k}F)e_k$ on smooth cylinder functions.
For finite signed measures $\sigma,\eta$ of finite Green energy, let
$U_D\sigma\in H$ be determined by
$\langle U_D\sigma,u\rangle_H=\int_Du\,d\sigma$ for $u\in H$ and write
\[
 I_G(\sigma,\eta)=\iint_{D\times D}G_D(x,z)\,\sigma(dx)\eta(dz)
 =\langle U_D\sigma,U_D\eta\rangle_H,
 \qquad I_G(\sigma)=I_G(\sigma,\sigma).
\]
Put $\nu_i=\phi_i^2\mu$, so that $\nu_i(D)=1$, and
$\Gamma F=\|\nabla_hF\|_H^2$.
By \cite[Lemma A.5]{CaiMV}, the Gaussian derivatives of the
logarithmic eigenvalues satisfy
\begin{equation}\label{eq:intro-response}
 \nabla_hY_i=-\gamma U_D\nu_i,\qquad
 \Gamma g=\gamma^2I_G(\nu_2-\nu_1).
\end{equation}
The Green-energy and spectral moment bounds from
\cite[Lemma A.1 and Proposition A.2]{CaiMV} give
\begin{equation}\label{eq:response-moments}
 \E[I_G(\nu_i)^s]<\infty,\qquad \E|Y_i|^s<\infty,
 \qquad s>0,\quad i\le N.
\end{equation}
Thus, for $u\in C_c^\infty(\R^N)|_S$, the Gaussian chain rule yields
\begin{equation}\label{eq:core-chain-rule}
 E_0(u,u)=\gamma^2\E\sum_{i,j\le N}
 \partial_i u(Y)\partial_ju(Y)I_G(\nu_i,\nu_j).
\end{equation}

For $0<\tau\le1/2$, define the heat operators and kernels
\begin{equation}\label{eq:heat-notation}
 Q_\tau e_k=e^{-\tau\kappa_k}e_k,\quad
 G_\tau(x,z)=\sum_{k\ge1}e^{-\tau\kappa_k}e_k(x)e_k(z),\quad
 \Gamma_\tau F=\langle Q_\tau\nabla_hF,\nabla_hF\rangle_H.
\end{equation}
The operators $Q_\tau$ differ from the weights $P_s$ in
\eqref{eq:H-basis}. If $p_t^D$ is the Dirichlet heat kernel, then
\begin{equation}\label{eq:heat-kernels}
 G_\tau(x,z)=2\pi\int_\tau^\infty p_t^D(x,z)\,dt,
 \qquad G_D-G_\tau=2\pi\int_0^\tau p_t^D\,dt.
\end{equation}
Both $G_\tau$ and $G_D-G_\tau$ are positive-definite and pointwise
nonnegative.

We use Kahane's comparison inequality and the positive and negative
moments of subcritical GMC from
\cite[Theorems 2.1, 2.11 and 2.12]{RV}. Put
\begin{equation}\label{eq:q-star}
 q_*=\frac4{\gamma^2},\qquad
 \zeta(q)=\left(2+\frac{\gamma^2}{2}\right)q
             -\frac{\gamma^2}{2}q^2.
\end{equation}
The dyadic estimate in \cite[proof of Lemma A.1]{CaiMV}, followed
by a covering with squares of side comparable to $r$, gives
\begin{equation}\label{eq:ball-moment}
 \sup_{x\in\overline D}\E\mu(B(x,r)\cap D)^q
 \le C_qr^{\zeta(q)},\qquad 1<q<q_*,\quad 0<r\le1.
\end{equation}
Since $\zeta(1)=2$ and $\zeta'(1)=2-\gamma^2/2>0$, fix
$1<a_0<q_*$ and $0<\alpha<(\zeta(a_0)-2)/a_0$ with
$\zeta(a_0)>2$. By \cite[Lemma A.1]{CaiMV}, there is a random
$C_\alpha\ge1+\mu(D)$ such that
\begin{equation}\label{eq:frostman}
 \mu(B(x,r)\cap D)\le C_\alpha r^\alpha
 \quad(x\in\overline D,\ 0<r\le1),\qquad
 \E C_\alpha^{a_0}<\infty.
\end{equation}
For a square $Q\subset D$ of side $r$ at distance comparable to $r$
from $\partial D$, the rescaled GMC comparison in
Appendix~\ref{sec:gmc} gives
\begin{equation}\label{eq:square-moment}
 \E\mu(Q)^q\le C_qr^{2q},\qquad 1<q<q_*.
\end{equation}
The constants are uniform for fixed comparability bounds.

\section{Spatial localization of the first gap}\label{sec:whitney}

Use a fixed Whitney decomposition $W$ subdivided a fixed number
of times; write $r_Q$ and $c_Q$ for the side length and center of $Q$.
After subdivision,
\begin{equation}\label{eq:whitney-geometry}
 c_0r_Q\le\operatorname{dist}(Q,\partial D)
 \le C_0r_Q,\qquad c_0>20.
\end{equation}
Each ball $U_Q=B(c_Q,2r_Q)$ is compactly contained in $D$; the
$U_Q$ have uniformly bounded overlap and the number of squares at scale
$r_j\asymp2^{-j}$ is $O(r_j^{-1})$.
The first-moment identity $\E\mu(A)=|A|$ for deterministic
measurable $A\subset D$ implies
\[
 \E\mu\left(\bigcup_{Q\in W}\partial Q\right)=0,
\]
so the Whitney squares partition $\mu$-almost all of $D$.
Say that $Q$ and $R$ are
separated, and write $Q\sim R$, if
\begin{equation}\label{eq:separation}
 \operatorname{dist}(Q,R)\ge\max(r_Q,r_R).
\end{equation}
Write $d(x)=\operatorname{dist}(x,\partial D)$.  The Green kernel on $D$ satisfies
\cite[Theorem 6.23]{CZ}
\begin{equation}\label{eq:green-comparison}
 G_D(x,z)\asymp_D
 \log\left(1+\frac{d(x)d(z)}{|x-z|^2}\right),
 \qquad x\ne z.
\end{equation}
For $x\in Q$ and $z\in R$, we have $d(x)\asymp r_Q$,
$d(z)\asymp r_R$ and $|x-z|\le\operatorname{diam}D$.  Thus
\begin{equation}\label{eq:green-pairs}
 \inf_{x\in Q,z\in R}G_D(x,z)\ge c_Dr_Qr_R,
 \qquad
 \inf_{x\in Q,z\in R}G_D(x,z)\ge c_D
 \quad\text{if }Q\not\sim R.
\end{equation}
Indeed, $Q\not\sim R$ implies $r_Q\asymp r_R$ and
$|x-z|\lesssim r_Q$ on $Q\times R$.  On the other hand, if
$Q\sim R$, then $|x-z|\ge\max(r_Q,r_R)$ and
\[
 \sup_{x\in Q,z\in R}G_D(x,z)
 \le C_D\log\left(1+C_D\frac{r_Qr_R}{\max(r_Q,r_R)^2}\right)
 \le C_D.
\]
All comparison constants are uniform in $Q,R$.  Set
\[
 \Lambda=1+\log C_\alpha,\qquad m_Q=\mu(Q).
\]

\begin{lemma}\label{lem:square-trace}
For every $W'\subset W$ and $f\in H$,
\begin{equation}\label{eq:trace}
 \sum_{Q\in W'}\int_Qf^2\,d\mu
 \le C\biggl(\sup_{Q\in W'}
 m_Q\bigl[\Lambda+\log_+(1/m_Q)\bigr]\biggr)\|f\|_H^2.
\end{equation}
The factor on the right is interpreted as zero when $m_Q=0$.
\end{lemma}

\begin{proof}
Fix $Q\in W$ with $m_Q>0$.  The Green kernel of $U_Q$ is
bounded on $Q\times Q$ by $C+\log(Cr_Q/|x-z|)$.  The Schur
constant of $\mu|_Q$ satisfies
\begin{align}
 V_Q&=\sup_{x\in Q}\int_QG_{U_Q}(x,z)\,\mu(dz)\notag\\
 &\le Cm_Q+C\int_0^{Cr_Q}
 \min\{m_Q,C_\alpha s^\alpha\}\,\frac{ds}{s}\notag\\
 &\le Cm_Q\left[1+\log_+\left(\frac{C_\alpha(Cr_Q)^\alpha}{m_Q}\right)\right]\notag\\
 &\le Cm_Q\bigl[\Lambda+\log_+(1/m_Q)\bigr].\label{eq:schur}
\end{align}
Choose $\eta_Q\in C_c^\infty(U_Q)$, equal to one on $Q$, with
$|\nabla\eta_Q|\le C/r_Q$.  The Green duality and Schur test give
\[
 \int_Q f^2\,d\mu\le
 V_Q\|\eta_Qf\|_{H_0^1(U_Q)}^2.
\]
Finite overlap and the Hardy inequality on $D$ yield
\begin{align*}
 \sum_{Q\in W'}\|\eta_Qf\|_{H_0^1(U_Q)}^2
 &\le C\int_D\left(|\nabla f|^2+\frac{f^2}{d(x)^2}\right)dx\\
 &\le C_D\|f\|_H^2.
\end{align*}
Taking the supremum of the Schur constants proves the assertion.
\end{proof}

The following elementary logarithmic Green estimate supplies the
spectral-window bound needed for the normalized nodal functions.

\begin{lemma}\label{lem:green-square}
For each $\ell>0$, on the event $\lambda_1\ge\ell$,
\begin{equation}\label{eq:green-square}
 \sup_{x\in D}\int_DG_D(x,z)^2\,\mu(dz)
 \le C_{\ell,D,\alpha}\,[1+\log C_\alpha].
\end{equation}
\end{lemma}

\begin{proof}
Fix $x\in D$ and put $v_R(z)=\min\{G_D(x,z),R\}$ for $R\ge1$.
With $\mathcal H^1$ denoting arclength, the Green flux identity
and the coarea formula give
\[
 \int_{\{G_D(x,\cdot)=t\}}|\nabla_zG_D(x,z)|\,d\mathcal H^1(z)
 =2\pi\quad\text{for a.e. }t>0,
 \qquad
 \|v_R\|_H^2=\frac1{2\pi}\int_0^R2\pi\,dt=R.
\]
The Rayleigh principle therefore gives
\[
 \int_Dv_R^2\,d\mu\le\lambda_1^{-1}\|v_R\|_H^2
 \le R/\ell.
\]
The kernel upper bound $G_D(x,z)\le C_D+\log^+(1/|x-z|)$ implies
\[
 \mu\{z:G_D(x,z)>t\}
 \le CC_\alpha e^{-\alpha t},\qquad t\ge C_D+1.
\]
By layer-cake integration,
\begin{align*}
 &\int_{\{G_D(x,\cdot)>R\}}G_D(x,z)^2\,\mu(dz)\\
 &\quad=R^2\mu\{G_D(x,\cdot)>R\}
   +2\int_R^\infty t\mu\{G_D(x,\cdot)>t\}\,dt\\
 &\quad\le C C_\alpha e^{-\alpha R}(1+R^2).
\end{align*}
Choose $R=C_D+1+(2/\alpha)\log(e+C_\alpha)$.  The displayed
estimates imply \eqref{eq:green-square} uniformly in $x$.
\end{proof}

Fix $0<\ell<L_2<\infty$ for the remainder of this section.  Define
\begin{equation}\label{eq:nodal-functions}
 f_+=\frac{(\phi_2)_+}{\|(\phi_2)_+\|_{L^2(\mu)}},\qquad
 f_-=\frac{(\phi_2)_-}{\|(\phi_2)_-\|_{L^2(\mu)}},\qquad
 p_Q^\pm=\int_Qf_\pm^2\,d\mu.
\end{equation}
Both denominators are nonzero because $\phi_2$ is orthogonal to the
strictly positive first eigenfunction.  Testing \eqref{eq:eigenpairs}
with the positive and negative parts shows
\begin{equation}\label{eq:nodal-energies}
 \|f_\pm\|_{L^2(\mu)}=1,\qquad
 \|f_\pm\|_H^2=\lambda_2,
 \qquad\langle f_+,f_-\rangle_{L^2(\mu)}=0.
\end{equation}
Write $K$ for the Green operator on $L^2(\mu)$, whose eigenvalues are
$\lambda_j^{-1}$.  The spectral Cauchy--Schwarz inequality and the Rayleigh principle
imply
\[
 \langle f_\pm,Kf_\pm\rangle
 \ge\frac{\|f_\pm\|_{L^2(\mu)}^4}{\|f_\pm\|_H^2}
 =\lambda_2^{-1},\qquad
 \lambda_1^{-1}\ge
 \left\langle\frac{f_++f_-}{\sqrt2},
 K\frac{f_++f_-}{\sqrt2}\right\rangle.
\]
It follows that
\begin{equation}\label{eq:nodal-coupling}
 \lambda_1^{-1}\ge\lambda_2^{-1}+\langle f_+,Kf_-\rangle,
 \qquad
 g\ge 1-\frac{\lambda_1}{\lambda_2}
    \ge\lambda_1\langle f_+,Kf_-\rangle.
\end{equation}
For $v\in C_c^\infty(D)$ with $v\ge0$, choose smooth
nondecreasing $\eta_\varepsilon$ with $\eta_\varepsilon=0$ on
$(-\infty,0]$ and $\eta_\varepsilon=1$ on
$[\varepsilon,\infty)$.  Testing \eqref{eq:eigenpairs} with
$\eta_\varepsilon(\phi_2)v$ gives
\begin{align*}
 \lambda_2\int_D\phi_2\eta_\varepsilon(\phi_2)v\,d\mu
 &=\frac1{2\pi}\int_D\eta_\varepsilon(\phi_2)
               \nabla\phi_2\cdot\nabla v\,dx\\
 &\quad+\frac1{2\pi}\int_D\eta_\varepsilon'(\phi_2)
                v|\nabla\phi_2|^2\,dx
 \ge\frac1{2\pi}\int_D\eta_\varepsilon(\phi_2)
                \nabla\phi_2\cdot\nabla v\,dx.
\end{align*}
Since
\[
 \eta_\varepsilon(\phi_2)\nabla\phi_2
   \longrightarrow\nabla(\phi_2)_+\quad\text{in }L^2(D),
 \qquad
 \phi_2\eta_\varepsilon(\phi_2)
   \longrightarrow(\phi_2)_+\quad\text{in }L^2(\mu),
\]
we may let $\varepsilon\downarrow0$ and repeat the argument with
$-\phi_2$.  The Rayleigh bound
$\|v\|_{L^2(\mu)}^2\le\lambda_1^{-1}\|v\|_H^2$
extends the resulting inequality from nonnegative
$C_c^\infty(D)$ functions to $H_+$:
\[
 \langle(\phi_2)_\pm,v\rangle_H
 \le\lambda_2\int_D(\phi_2)_\pm v\,d\mu,
 \qquad v\in H,\quad v\ge0.
\]
For $w_\pm=\lambda_2K(\phi_2)_\pm-(\phi_2)_\pm$, Green
duality gives $\langle w_\pm,v\rangle_H\ge0$ for $v\ge0$.
Taking $v=(w_\pm)_-$ gives
\[
 0\le\langle w_\pm,(w_\pm)_-\rangle_H
 =-\|(w_\pm)_-\|_H^2,
 \qquad f_\pm\le\lambda_2Kf_\pm.
\]
By Cauchy--Schwarz and Lemma~\ref{lem:green-square}, on
$\{\lambda_1\ge\ell,\lambda_2\le L_2\}$,
\begin{equation}\label{eq:nodal-sup}
 \|f_+\|_{L^\infty(\mu)}^2+\|f_-\|_{L^\infty(\mu)}^2
 \le C_{\ell,L_2}(1+\log C_\alpha).
\end{equation}

Fix $\varepsilon\in(0,e^{-1})$ sufficiently small, depending only on
$D,\ell,L_2,\alpha$, and put
\begin{equation}\label{eq:heavy-squares}
 W_{\mathrm{heavy}}=\{Q\in W:m_Q\ge\varepsilon\Lambda^{-2}\}.
\end{equation}
Since $\Lambda\ge1$, for $0<t<\varepsilon\Lambda^{-2}$,
\[
 t\bigl[\Lambda+\log_+(1/t)\bigr]
 \le C\varepsilon\bigl[1+\log(1/\varepsilon)\bigr].
\]
Lemma~\ref{lem:square-trace} and \eqref{eq:nodal-energies} yield
\[
 \sum_{Q\notin W_{\mathrm{heavy}}}p_Q^\pm
 \le C L_2\varepsilon[1+\log(1/\varepsilon)]\le\frac18.
\]
Hence
\begin{equation}\label{eq:heavy-mass}
 \sum_{Q\in W_{\mathrm{heavy}}}p_Q^+\ge\frac78,
 \qquad\sum_{Q\in W_{\mathrm{heavy}}}p_Q^-\ge\frac78.
\end{equation}

\begin{lemma}\label{lem:gap-geometry}
For $p>0$, on $\{\lambda_1\ge\ell,\lambda_2\le L_2\}$,
\begin{equation}\label{eq:gap-geometry}
 g^{-p}\le C_{\ell,L_2,p}\Lambda^p
 \left[1+\sum_{\substack{Q\sim R\\Q,R\in W_{\mathrm{heavy}}}}
 (r_Qr_R)^{-p}\right].
\end{equation}
\end{lemma}

\begin{proof}
Put $M^2=C_{\ell,L_2}\Lambda$.  By \eqref{eq:nodal-sup},
\[
 \int_Q f_\pm\,d\mu\ge M^{-1}p_Q^\pm.
\]
Using the Whitney partition of $\mu$ and
\eqref{eq:nodal-coupling},
\[
 \langle f_+,Kf_-\rangle
 =\sum_{Q,R\in W}\iint_{Q\times R}
 G_D(x,z)f_+(x)f_-(z)\,\mu(dx)\mu(dz).
\]
Thus \eqref{eq:green-pairs} gives
\begin{equation}\label{eq:coupling-pairs}
 g\ge\frac{c\ell}{M^2}
 \sum_{Q,R\in W}\inf_{Q\times R}G_D\,p_Q^+p_R^-.
\end{equation}
In particular,
\begin{equation}\label{eq:nonseparated-weight}
 \sum_{\substack{Q,R\in W_{\mathrm{heavy}}\\Q\not\sim R}}
 p_Q^+p_R^-\le CgM^2/\ell.
\end{equation}
If $CgM^2/\ell>1/4$, then $g^{-p}\le CM^{2p}$.
Otherwise \eqref{eq:heavy-mass} yields
\[
 s_0:=\sum_{\substack{Q,R\in W_{\mathrm{heavy}}\\Q\sim R}}
 p_Q^+p_R^-\ge\frac{49}{64}-\frac14>\frac12.
\]
Set $w_{QR}=p_Q^+p_R^-/s_0$ on these pairs.  Then
\begin{align*}
 g&\ge\frac{c\ell s_0}{M^2}\sum_{Q\sim R}w_{QR}r_Qr_R,\\
 g^{-p}&\le CM^{2p}\left(\sum_{Q\sim R}w_{QR}r_Qr_R\right)^{-p}\\
 &\le CM^{2p}\sum_{Q\sim R}w_{QR}(r_Qr_R)^{-p}\\
 &\le C\Lambda^p\sum_{\substack{Q\sim R\\Q,R\in W_{\mathrm{heavy}}}}
 (r_Qr_R)^{-p}.
\end{align*}
The middle inequality is Jensen's inequality.
\end{proof}

\section{Gaussian rotations and inverse-gap moments}\label{sec:rotations}

Fix $0<\tau\le1/2$.  Let $h_\tau^{\mathrm{hi}},U_\tau,V_\tau$ be
independent centered Gaussian fields of covariances $G_D-G_\tau$,
$G_\tau$ and $G_\tau$, respectively.  Put
\begin{equation}\label{eq:rotation}
 h_\theta=h_\tau^{\mathrm{hi}}+U_\tau\cos\theta+V_\tau\sin\theta,
 \quad
 \dot h_\theta=-U_\tau\sin\theta+V_\tau\cos\theta,
 \quad 0\le\theta\le2\pi.
\end{equation}
Expectations on this auxiliary product space are also denoted by $\E$.
The rotated GMC measure $\mu_\theta=M_\gamma(h_\theta)$ is defined
using the covariance $G_D$.  Equivalently,
\begin{equation}\label{eq:rotated-chaos}
 \mu_\theta(dx)=
 \exp\left(\gamma[U_\tau(x)\cos\theta+V_\tau(x)\sin\theta]
 -\frac{\gamma^2}{2}G_\tau(x,x)\right)
 M_\gamma(h_\tau^{\mathrm{hi}})(dx).
\end{equation}
For each fixed $\theta$, $h_\theta$ has the original GFF law and
$\dot h_\theta$ is independent of $h_\theta$ with covariance $G_\tau$.
The low-frequency fields are almost surely smooth and bounded.
For bounded measurable $f$,
\begin{equation}\label{eq:rotated-derivative}
 \frac d{d\theta}\int_D f\,d\mu_\theta
 =\gamma\int_D f\dot h_\theta\,d\mu_\theta.
\end{equation}
For $\theta,\varphi\in[0,2\pi]$, put
\[
 c_{\theta,\varphi}
 =\gamma\bigl\|U_\tau(\cos\theta-\cos\varphi)
       +V_\tau(\sin\theta-\sin\varphi)\bigr\|_{L^\infty(D)}.
\]
Almost surely, $c_{\theta,\varphi}\to0$ as
$\theta\to\varphi$, and
\[
 e^{-c_{\theta,\varphi}}\mu_\varphi
 \le\mu_\theta\le
 e^{c_{\theta,\varphi}}\mu_\varphi.
\]
The min--max principle therefore gives
\[
 e^{-c_{\theta,\varphi}}\lambda_i(\mu_\varphi)
 \le\lambda_i(\mu_\theta)\le
 e^{c_{\theta,\varphi}}\lambda_i(\mu_\varphi),
 \qquad i\le N.
\]
Thus the spectral coordinates are continuous along the circle.

For $Q\in W$, define
\begin{equation}\label{eq:scaled-mass}
 M_Q=r_Q^{-2}\mu(Q),\qquad
 A_Q=\frac{I_G(\mu|_Q)}{\mu(Q)^2},
\end{equation}
with $A_Q=0$ when $\mu(Q)=0$.

\begin{lemma}\label{lem:rotation-pair}
For every $1<q<q_*$ and $Q\sim R$,
\begin{equation}\label{eq:rotation-pair}
 \sup_{0<\tau\le1/2}
 \E\left[\sup_{\theta\in[0,2\pi]}
 (\mu_\theta(Q)\mu_\theta(R))^q\right]
 \le C_q(r_Qr_R)^{2q}.
\end{equation}
\end{lemma}

\begin{lemma}\label{lem:rotation-frostman}
For $a_0,\alpha$ fixed in \eqref{eq:frostman}, there are random
variables $C_\tau\ge1+\sup_\theta\mu_\theta(D)$ such that
\begin{equation}\label{eq:rotation-frostman}
 \mu_\theta(B(x,r)\cap D)\le C_\tau r^\alpha
 \quad(\theta\in[0,2\pi],\ x\in\overline D,\ 0<r\le1),
 \qquad
 \sup_{0<\tau\le1/2}\E C_\tau^{a_0}<\infty.
\end{equation}
\end{lemma}

The two uniform estimates are proved in Appendix~\ref{sec:rotation-technical}.

\subsection{The inverse-gap moment}

Fix $0<\ell<L_2<\infty$.  With the common Frostman constant
$C_\tau$ from Lemma~\ref{lem:rotation-frostman}, put
\[
 \Lambda_\tau=1+\log C_\tau.
\]
At each angle define $W_\theta$ by
\eqref{eq:heavy-squares}, with $C_\alpha=C_\tau$.
Define
\begin{equation}\label{eq:circle-inverse-gap}
 R_\tau=\max\left\{1,
 \sup_{\substack{0\le\theta\le2\pi:\ \lambda_1(\mu_\theta)\ge\ell\\
 \lambda_2(\mu_\theta)\le L_2}}
 \frac1{\log(\lambda_2(\mu_\theta)/\lambda_1(\mu_\theta))}
 \right\}.
\end{equation}
The supremum over the empty set is zero.  Almost surely the
measures $\mu_\theta$ are mutually equivalent through bounded positive
densities.  Positivity improvement of the Green operator on connected
$D$ gives $\lambda_1(\mu_\theta)<\lambda_2(\mu_\theta)$
for all $\theta$.  Continuity of the eigenvalues along the compact
circle then implies $R_\tau<\infty$.

\begin{proposition}\label{prop:circle-negative-moments}
For every $0<p<2q_*-1=8/\gamma^2-1$,
\begin{equation}\label{eq:circle-negative-moments}
 \sup_{0<\tau\le1/2}\E R_\tau^p<\infty.
\end{equation}
\end{proposition}

\begin{proof}
At each angle, a pair $Q,R\in W_\theta$ satisfies
$\mu_\theta(Q)\mu_\theta(R)\ge
\varepsilon^2\Lambda_\tau^{-4}$.  Hence, for $q>0$,
\[
 \one_{\{Q,R\in W_\theta\}}
 \le\varepsilon^{-2q}\Lambda_\tau^{4q}
       [\mu_\theta(Q)\mu_\theta(R)]^q.
\]
Lemma~\ref{lem:gap-geometry}, followed by the supremum over the
admissible angles, yields
\begin{equation}\label{eq:circle-gap-bound}
 R_\tau^p
 \le C\Lambda_\tau^p+C\Lambda_\tau^{p+4q}
 \sum_{Q\sim R}(r_Qr_R)^{-p}
 \sup_\theta[\mu_\theta(Q)\mu_\theta(R)]^q.
\end{equation}
Choose
\[
 \max\{1,(p+1)/2\}<q<q'<q_*.
\]
For $q<q'<q_*$, H\"older's inequality gives
\begin{align*}
 &\E\left[\Lambda_\tau^{p+4q}
       \sup_\theta(\mu_\theta(Q)\mu_\theta(R))^q\right]\\
 &\quad\le
 \left(\E\Lambda_\tau^{(p+4q)q'/(q'-q)}\right)^{1-q/q'}
 \left(\E\sup_\theta
       (\mu_\theta(Q)\mu_\theta(R))^{q'}\right)^{q/q'}
 \le C(r_Qr_R)^{2q}.
\end{align*}
The two factors are controlled by Lemmas~\ref{lem:rotation-frostman}
and \ref{lem:rotation-pair}.  Summing by Whitney scale gives
\begin{align*}
 \E R_\tau^p
 &\le C+C\sum_{j,k\ge0}r_j^{-1}r_k^{-1}
 (r_jr_k)^{2q-p}\\
 &\le C+C\left(\sum_{j\ge0}r_j^{2q-p-1}\right)^2<\infty.
\end{align*}
\end{proof}

\section{Response energy and capacity}\label{sec:curvature}

For $n=1,2$, write $J_n=\sum_{i=1}^nY_i$.
Fix $0<\tau\le1/2$ and a realization of $h_\tau^{\mathrm{hi}}$.
Set
\[
 \mu_0(dx)=e^{-\gamma^2G_\tau(x,x)/2}
                       M_\gamma(h_\tau^{\mathrm{hi}})(dx),
 \qquad \mu_w=e^{\gamma w}\mu_0
\]
for bounded $w$, and let $J_n(w)$ be the logarithmic spectral sum
of $\mu_w$.  For an $n$-dimensional subspace $V\subset H$ with
$H$-orthonormal basis $v_1,\ldots,v_n$, the min--max principle gives
\begin{equation}\label{eq:logdet-concavity}
 J_n(w)=\inf_{\dim V=n}
 \left[-\log\det\left(\int_Dv_av_b e^{\gamma w}\,d\mu_0
 \right)_{a,b\le n}\right],
\end{equation}
The Gram determinant admits the representation
\begin{equation}\label{eq:alternant}
 \det\!\left(\int_Dv_av_b e^{\gamma w}\,d\mu_0\right)_{a,b\le n}
 =\frac1{n!}\int_{D^n}
 \det\bigl(v_a(x_b)\bigr)_{a,b\le n}^{\,2}
 e^{\gamma\sum_{b=1}^n w(x_b)}\prod_{b=1}^n\mu_0(dx_b).
\end{equation}
H\"older's inequality makes the logarithm of the integral in
\eqref{eq:alternant} convex in $w$.  Taking the infimum of the
corresponding concave functions yields, for $0\le t\le1$,
\[
 J_n(tw+(1-t)\widetilde w)
 \ge tJ_n(w)+(1-t)J_n(\widetilde w).
\]
Furthermore,
\begin{equation}\label{eq:sn-lipschitz}
 |J_n(w)-J_n(\widetilde w)|
 \le n\gamma\|w-\widetilde w\|_\infty.
\end{equation}

Along the rotation \eqref{eq:rotation}, let $J_n(\theta)$ denote the
corresponding spectral function and put
\begin{equation}\label{eq:field-sup}
 T_\tau=
 \|U_\tau\|_\infty+
 \|V_\tau\|_\infty.
\end{equation}
The fields $U_\tau,V_\tau$ are those in
\eqref{eq:rotation}.  Put
$w_\theta=U_\tau\cos\theta+V_\tau\sin\theta$.
Since
\[
 \frac{w_{\theta+u}+w_{\theta-u}}2=(\cos u)w_\theta,
\]
concavity and \eqref{eq:sn-lipschitz} yield
\begin{align}
 &J_n(\theta+u)+J_n(\theta-u)-2J_n(\theta)\notag\\
 &\quad\le2[J_n((\cos u)w_\theta)-J_n(w_\theta)]\notag\\
 &\quad\le2n\gamma(1-\cos u)\|w_\theta\|_\infty
 \le n\gamma T_\tau u^2.
 \label{eq:semiconcavity}
\end{align}
Thus the distributional second derivative is a finite signed measure
satisfying $J_n''\le n\gamma T_\tau\,d\theta$.
Periodicity gives $J_n''(\mathbb T)=0$.  Since $g=J_2-2J_1$,
\begin{equation}\label{eq:sn-total-variation}
 \|J_n''\|_{\TV}\le4\pi n\gamma T_\tau,\qquad
 \|g''\|_{\TV}\le16\pi\gamma T_\tau.
\end{equation}
For $s>0$, we have
\begin{equation}\label{eq:field-sup-moment}
 \|T_\tau\|_{L^s(\Pp)}\le C_s[1+\log(1/\tau)].
\end{equation}
To verify this, put $L_\tau=1+\log(1/\tau)$.  The heat-kernel
bounds give, for a fixed $J_0<\infty$,
\[
 \sup_x\E U_\tau(x)^2\le CL_\tau,\qquad
 \E|U_\tau(x)-U_\tau(z)|^2
 \le C\tau^{-J_0}|x-z|^2.
\]
Fix $a>J_0/2+1$.  Partition a square containing $\overline D$
into nested dyadic cells of side $2^{-j}\tau^a$, and select a
point of $\overline D$ from each cell meeting $\overline D$.
The resulting nets $\mathcal N_j$ satisfy
$|\mathcal N_j|\le C\tau^{-2a}2^{2j}$.
The parent cells define maps $z_{j-1}:\mathcal N_j\to\mathcal N_{j-1}$
with $|x-z_{j-1}(x)|\le C2^{-j}\tau^a$.  For every
$x\in\overline D$ there is a chain $x_j\in\mathcal N_j$ such that
$x_j\to x$ and $x_{j-1}=z_{j-1}(x_j)$.
For centered Gaussian variables $Z_1,\ldots,Z_M$ and $s\ge1$,
\[
 \left\|\max_{i\le M}|Z_i|\right\|_{L^s}
 \le C_s\sqrt{1+\log M}\max_{i\le M}\|Z_i\|_{L^2}.
\]
Consequently,
\begin{align*}
 \left\|\max_{x\in\mathcal N_0}|U_\tau(x)|\right\|_{L^s}
 &\le C_s L_\tau,\\
 \left\|\max_{x\in\mathcal N_j}
   |U_\tau(x)-U_\tau(z_{j-1}(x))|\right\|_{L^s}
 &\le C_s\tau^{a-J_0/2}2^{-j}\sqrt{L_\tau+j},\quad j\ge1.
\end{align*}
Continuity and the identity
\[
 U_\tau(x)=U_\tau(x_0)
 +\sum_{j\ge1}[U_\tau(x_j)-U_\tau(x_{j-1})]
\]
then give
\[
 \|\|U_\tau\|_\infty\|_{L^s}
 \le C_s\left[L_\tau+
 \tau^{a-J_0/2}\sum_{j\ge1}2^{-j}\sqrt{L_\tau+j}\right]
 \le C_sL_\tau.
\]
The same estimate holds for $V_\tau$; the case $0<s<1$ follows
from the $L^1$ bound.

Fix a nonzero $\chi\in C_c^\infty(\R)$ with $0\le\chi\le1$.
Choose $\omega\in C_c^\infty(\R)$ with $0\le\omega\le1$,
$\omega=1$ on $[-1,1]$ and $\omega=0$ outside $[-2,2]$.
Define the spectral cutoff $\psi(y)=\chi(y_1)\omega(y_2-y_1)$
and along the rotation put
\begin{equation}\label{eq:Psi}
 \Psi(\theta)=\psi(Y(h_\theta)).
\end{equation}
Set
\[
 \ell=\min_{y\in\supp\chi}e^y,
 \qquad L_2=e^2\max_{y\in\supp\chi}e^y.
\]
On $\{\Psi\ne0\}$ we have $\lambda_1\ge\ell$ and
$\lambda_2\le L_2$.  Use the corresponding $R_\tau$ from
\eqref{eq:circle-inverse-gap}.

\begin{lemma}\label{lem:smoothed-energy}
For every such $\chi$,
\begin{equation}\label{eq:smoothed-energy}
 \E\left[\chi(Y_1)^2\frac{\Gamma_\tau g}{g^2};\ 0<g<1\right]
 \le C_\chi[1+\log(1/\tau)],
 \qquad 0<\tau\le\frac12.
\end{equation}
\end{lemma}

\begin{proof}
Fix one sample of the rotation.  Since $g(\theta)>0$ on the compact
circle, integration by parts for the signed measure $g''$ gives
\begin{align}
 \int_0^{2\pi}\Psi^2\frac{(g')^2}{g^2}\,d\theta
 &=\int_0^{2\pi}\frac{\Psi^2}{g}\,dg'
   +2\int_0^{2\pi}\frac{\Psi\Psi'g'}g\,d\theta\notag\\
 &\le R_\tau\|g''\|_{\TV}
   +\frac12\int_0^{2\pi}\Psi^2\frac{(g')^2}{g^2}\,d\theta
   +2\int_0^{2\pi}|\Psi'|^2\,d\theta.
\end{align}
By the fixed-angle simplicity in \cite{Cai} and Fubini's
theorem, the first-order response applies for almost every angle.
The integration by parts above holds because $g>0$ is Lipschitz
and $g'$ has bounded variation.
Consequently,
\begin{equation}\label{eq:circle-energy}
 \int_0^{2\pi}\Psi^2\frac{(g')^2}{g^2}\,d\theta
 \le2R_\tau\|g''\|_{\TV}
 +4\int_0^{2\pi}|\Psi'|^2\,d\theta.
\end{equation}
Conditionally on $h_\theta$, the derivative $g'(\theta)$ is
$-\gamma\int_D\dot h_\theta\,d(\nu_2-\nu_1)$, whence
\begin{equation}\label{eq:rotation-response}
 \E[|g'(\theta)|^2\mid h_\theta]
 =\Gamma_\tau g(h_\theta).
\end{equation}
The same formula applies to the spectral cutoff $\psi$ and its
energy $\Gamma_\tau(\psi(Y))$.
Taking expectations in \eqref{eq:circle-energy}, integrating over
$\theta$, and using stationarity of the marginal $h_\theta$ gives
\begin{align*}
 2\pi\E\left[\psi(Y)^2\frac{\Gamma_\tau g}{g^2}\right]
 &\le 2\E[R_\tau\|g''\|_{\TV}]
      +8\pi\E\Gamma_\tau(\psi(Y))\\
 &\le C\E[R_\tau T_\tau]
      +8\pi\E\Gamma(\psi(Y)).
\end{align*}
Choose $p>1$ permitted by Proposition~\ref{prop:circle-negative-moments}.
H\"older's inequality and \eqref{eq:field-sup-moment} imply
\[
 \E[R_\tau T_\tau]
 \le C[1+\log(1/\tau)].
\]
The quantity $\E\Gamma(\psi(Y))$ is finite by \eqref{eq:response-moments}
and smoothness of the cutoffs.  Since $\omega=1$ for $0<g<1$, the left
side dominates the expression in \eqref{eq:smoothed-energy}.
\end{proof}

\subsection{The high-frequency remainder}

The preceding estimate concerns $\Gamma_\tau$ rather than
$\Gamma$.  By \eqref{eq:heat-kernels}, the Gaussian
heat-kernel bound and \eqref{eq:frostman} imply
\begin{equation}\label{eq:small-time-mass}
 \sup_x\int_D(G_D-G_\tau)(x,z)\,\mu(dz)
 \le CC_\alpha\tau^{\alpha/2}.
\end{equation}
Indeed, the Gaussian heat-kernel upper bound and dyadic annuli give
\[
 \begin{aligned}
 \int_Dp_t^D(x,z)\,\mu(dz)
 &\le Ct^{-1}\sum_{j\ge0}e^{-c4^j}
   \mu(B(x,2^{j+1}\sqrt t)\cap D)\\
 &\le CC_\alpha t^{\alpha/2-1}
   \sum_{j\ge0}2^{j\alpha}e^{-c4^j}
 \le CC_\alpha t^{\alpha/2-1}.
 \end{aligned}
\]
Integration over $0<t<\tau$ proves \eqref{eq:small-time-mass}.

\begin{lemma}\label{lem:high-frequency}
For each $0<\ell<L_2<\infty$ there exist $C,\delta>0$ such that
\begin{equation}\label{eq:high-frequency}
 \E\left[\one_{\{\lambda_1\ge\ell,\,\lambda_2\le L_2\}}
 (\Gamma g-\Gamma_\tau g)\right]
 \le C\tau^\delta,\qquad0<\tau\le\frac12.
\end{equation}
\end{lemma}

\begin{proof}
The difference kernel is pointwise nonnegative and positive definite.
Writing $\sigma=\nu_2-\nu_1$ and using $\nu_i\ge0$,
\begin{align*}
 0\le\Gamma g-\Gamma_\tau g
 &=\gamma^2I_{G_D-G_\tau}(\sigma)\\
 &\le\gamma^2I_{G_D-G_\tau}(\nu_1+\nu_2).
\end{align*}
On the indicated spectral window, Cauchy--Schwarz and
Lemma~\ref{lem:green-square} give
\[
 \sum_{i=1}^2\|\phi_i\|_{L^\infty(\mu)}^2
 \le C(1+\log C_\alpha).
\]
Since $\nu_1(D)+\nu_2(D)=2$,
\[
 I_{G_D-G_\tau}(\nu_1+\nu_2)
 \le2\left(\sum_{i=1}^2\|\phi_i\|_{L^\infty(\mu)}^2\right)
 \sup_x\int_D(G_D-G_\tau)(x,z)\,\mu(dz).
\]
Thus, by \eqref{eq:small-time-mass},
\begin{equation}\label{eq:high-frequency-pointwise}
 0\le\Gamma g-\Gamma_\tau g
 \le C(1+\log C_\alpha)C_\alpha\tau^{\alpha/2}.
\end{equation}
Put $\Delta_\tau=\Gamma g-\Gamma_\tau g$.  On the
spectral window in \eqref{eq:high-frequency}, the bounds
$0\le\Delta_\tau\le\Gamma g$ and
\eqref{eq:high-frequency-pointwise} imply
\[
 \Delta_\tau\le (\Gamma g)^{1-\vartheta}
 [C(1+\log C_\alpha)C_\alpha]^\vartheta
 \tau^{\alpha\vartheta/2}.
\]
Choose $0<\vartheta<\min\{1,a_0/2\}$ and
$0<\eta<a_0-2\vartheta$.  H\"older's inequality gives
\begin{align*}
 &\E[\one_{\{\lambda_1\ge\ell,\lambda_2\le L_2\}}\Delta_\tau]\\
 &\quad\le C\tau^{\alpha\vartheta/2}
 \bigl(\E(\Gamma g)^{2(1-\vartheta)}\bigr)^{1/2}
 \bigl(\E[C_\alpha^{2\vartheta}
       (1+\log C_\alpha)^{2\vartheta}]\bigr)^{1/2}
 \le C\tau^{\alpha\vartheta/2},
\end{align*}
since $C_\alpha^{2\vartheta}(1+\log C_\alpha)^{2\vartheta}
\le C_{\vartheta,\eta}C_\alpha^{2\vartheta+\eta}$ and
\eqref{eq:frostman}, \eqref{eq:response-moments} apply.
\end{proof}

\subsection{Spectral cutoff functions}\label{sec:capacity-cutoffs}

For $m\ge1$, choose $\chi_m\in C_c^\infty(\R)$ with
$0\le\chi_m\le1$ and $\chi_m=1$ on an open neighborhood of
$[-m,m]$.  For $n\ge1$, define
\begin{equation}\label{eq:capacity-cuts}
 v_n(r)=\begin{cases}
 1,&r\le e^{-n},\\
 -\log r/n,&e^{-n}<r<1,\\
 0,&r\ge1,
 \end{cases}
 \qquad
 u_{m,n}(y)=\chi_m(y_1)v_n(y_2-y_1).
\end{equation}

\begin{lemma}\label{lem:capacity-domain}
For every $m,n\ge1$,
\begin{equation}\label{eq:capacity-energy}
 u_{m,n}\in F_0^{\min},\qquad
 E_0^{\min}(u_{m,n},u_{m,n})
 =\E\Gamma(u_{m,n}(Y)).
\end{equation}
\end{lemma}

\begin{proof}
Fix $m,n$.  Choose $\eta_R\in C_c^\infty(\R^N)$ with
$\eta_R=1$ on $B_R$ and $|\nabla\eta_R|\le C/R$.
By the chain rule and \eqref{eq:response-moments},
\begin{align*}
 \|u_{m,n}\eta_R-u_{m,n}\|_{L^2(\rho)}^2&\longrightarrow0,\\
 \E\Gamma(((u_{m,n}\eta_R-u_{m,n})(Y)))
 &\le C_{m,n}\E\left[\one_{\{|Y|>R\}}
             \sum_{i\le N}\Gamma Y_i\right]
 +\frac C{R^2}\sum_{i\le N}\E\Gamma Y_i
 \longrightarrow0.
\end{align*}
For fixed $R$, let $u_{m,n,R,\varepsilon}$ be a smooth compactly
supported convolution approximation of $u_{m,n}\eta_R$.
The derivative discontinuities of $v_n$ occur only on
\[
 \{y:y_2-y_1=e^{-n}\}\cup\{y:y_2-y_1=1\}.
\]
These sets have $\rho$-measure zero by the joint spectral density \cite{Cai}.
The uniform Lipschitz bound and dominated convergence imply
\begin{align*}
 &\lim_{\varepsilon\downarrow0}
 \|u_{m,n,R,\varepsilon}-u_{m,n}\eta_R\|_{L^2(\rho)}^2=0,\\
 &\lim_{\varepsilon\downarrow0}
 \E\Gamma(((u_{m,n,R,\varepsilon}
                    -u_{m,n}\eta_R)(Y)))=0.
\end{align*}
The defining $C_c^\infty(\R^N)$ closure then gives
\eqref{eq:capacity-energy}.
\end{proof}

\subsection{Final proof of Theorem~\ref{thm:noncollision}}\label{sec:final-proof}

Fix $m\ge1$.  The functions $u_{m,n}$ belong to $F_0^{\min}$ by
Lemma~\ref{lem:capacity-domain}.  On $\{g>0\}$, the Gaussian chain
rule and \eqref{eq:capacity-cuts} give
\begin{align}
 \nabla_h(u_{m,n}(Y))
 &=\chi_m(Y_1)v_n'(g)\nabla_hg
   +v_n(g)\chi_m'(Y_1)\nabla_hY_1,\notag\\
 |v_n'(g)|^2
 &=\frac{\one_{\{e^{-n}<g<1\}}}{n^2g^2},
 \label{eq:final-cutoff-derivative}\\
 E_0^{\min}(u_{m,n},u_{m,n})
 &\le\frac2{n^2}
 \E\left[\chi_m(Y_1)^2\frac{\Gamma g}{g^2};\
                  e^{-n}<g<1\right]
 +2\E\left[v_n(g)^2\chi_m'(Y_1)^2\Gamma Y_1\right].
 \label{eq:capacity-main-energy}
\end{align}
Set
\[
 \ell_m=\min_{y\in\supp\chi_m}e^y>0,
 \qquad L_m=e\max_{y\in\supp\chi_m}e^y<\infty.
\]
In particular,
\begin{equation}\label{eq:final-spectral-window}
 \{\chi_m(Y_1)\ne0,\ 0<g<1\}
 \subset\{\lambda_1\ge\ell_m,\ \lambda_2\le L_m\}.
\end{equation}
Since $\Gamma_\tau g\le\Gamma g$, the decomposition
\[
 \Gamma g=\Gamma_\tau g+(\Gamma g-\Gamma_\tau g)
\]
and Lemmas~\ref{lem:smoothed-energy} and
\ref{lem:high-frequency}, applied with
\eqref{eq:final-spectral-window}, yield
\begin{align}
 &\E\left[\chi_m(Y_1)^2\frac{\Gamma g}{g^2};\
                       e^{-n}<g<1\right]\notag\\
 &\quad\le
 \E\left[\chi_m(Y_1)^2\frac{\Gamma_\tau g}{g^2};\
                       0<g<1\right]
 +e^{2n}\E\left[\one_{\{\lambda_1\ge\ell_m,\,\lambda_2\le L_m\}}
                  (\Gamma g-\Gamma_\tau g)\right]\notag\\
 &\quad\le C_m\bigl[1+\log(1/\tau)+e^{2n}\tau^{\delta_m}\bigr],
 \qquad 0<\tau\le\tfrac12,
 \label{eq:final-low-high-bound}
\end{align}
where $\delta_m>0$.  Substituting into
\eqref{eq:capacity-main-energy} gives
\begin{equation}\label{eq:capacity-bound}
 E_0^{\min}(u_{m,n},u_{m,n})
 \le\frac{C_m}{n^2}
       \bigl[1+\log(1/\tau)+e^{2n}\tau^{\delta_m}\bigr]
 +2\E\left[v_n(g)^2\chi_m'(Y_1)^2\Gamma Y_1\right].
\end{equation}
Choose $a_m>2/\delta_m$ and $\tau_n=e^{-a_m n}$ for all
sufficiently large $n$.  Then
\begin{align}
 \frac{1+\log(1/\tau_n)}{n^2}
 &=\frac{1+a_mn}{n^2}\longrightarrow0,
 \qquad
 \frac{e^{2n}\tau_n^{\delta_m}}{n^2}
 =\frac{e^{-(a_m\delta_m-2)n}}{n^2}\longrightarrow0.
 \label{eq:final-choice}
\end{align}
Static simplicity gives $g>0$ $\Pp$-almost surely, so
$v_n(g)\to0$ almost surely and $0\le v_n(g)\le1$.
By \eqref{eq:response-moments},
\begin{align}
 &\E\left[v_n(g)^2\chi_m'(Y_1)^2\Gamma Y_1\right]
 \longrightarrow0,\notag\\
 &\|u_{m,n}\|_{L^2(\rho)}^2
 =\E[\chi_m(Y_1)^2v_n(g)^2]\longrightarrow0.
 \label{eq:final-dominated}
\end{align}
Consequently, \eqref{eq:capacity-bound}--\eqref{eq:final-dominated}
prove the required convergence in the full form norm:
\begin{equation}\label{eq:capacity-convergence}
 \lim_{n\to\infty}
 \left[E_0^{\min}(u_{m,n},u_{m,n})
       +\|u_{m,n}\|_{L^2(\rho)}^2\right]=0.
\end{equation}

Let
\[
 A_m=\Sigma_1\cap\{y\in S:|y_1|\le m\}.
\]
As $\chi_m=1$ on a neighborhood of $[-m,m]$ and $v_n=1$
on $(-\infty,e^{-n}]$, for each $n$ there exists a relative
open neighborhood $O_{m,n}$ of $A_m$ in $S$ such that
\begin{equation}\label{eq:final-admissible}
 u_{m,n}(y)=1\quad(y\in O_{m,n}),\qquad
 A_m\subset O_{m,n}\subset S.
\end{equation}
The definition of $1$-capacity, Lemma~\ref{lem:capacity-domain},
and \eqref{eq:capacity-convergence} imply
\begin{equation}\label{eq:final-local-capacity}
 0\le\Capacity_{E_0^{\min},1}(A_m)
 \le\liminf_{n\to\infty}
       \bigl[E_0^{\min}(u_{m,n},u_{m,n})
             +\|u_{m,n}\|_{L^2(\rho)}^2\bigr]=0.
\end{equation}
Since $\Sigma_1=\bigcup_{m\ge1}A_m$, countable subadditivity
of the regular minimal-form capacity and \eqref{eq:final-local-capacity}
give $\Capacity_{E_0^{\min},1}(\Sigma_1)=0$.
By \eqref{eq:form-interface}, admissible functions for the minimal
form are also admissible for the canonical form, with the same energy.
Hence
\begin{equation}\label{eq:final-global-capacity}
 0\le\Capacity_{E_*,1}(\Sigma_1)
 \le\Capacity_{E_0^{\min},1}(\Sigma_1)=0.
\end{equation}

Let $\sigma_{\Sigma_1}$ denote the first hitting time of
$\Sigma_1$ by the conservative continuous Hunt diffusion
$X^{\min}$ of the regular minimal form.  The capacity--polarity theorem for regular symmetric Dirichlet forms
\cite[Chapter 4]{FOT} and
\eqref{eq:final-global-capacity} give
\[
 \Pp_x(\sigma_{\Sigma_1}<\infty)=0
 \quad\text{for $E_0^{\min}$-quasi-every }x\in S.
\]
The exceptional set has $\rho$-measure zero.  Since
$\rho(\Sigma_1)=0$, for every $\pi\ll\rho$,
\[
 \Pp_\pi(\sigma_{\Sigma_1}<\infty)=0,
 \qquad
 \Pp_\pi\{X^{\min}_1(0)<X^{\min}_2(0)\}=1.
\]
By path continuity, any loss of strict ordering forces a hit
of $\Sigma_1$.  Therefore
\begin{equation}\label{eq:final-noncollision}
 \Pp_\pi\{X^{\min}_1(t)<X^{\min}_2(t)
                    \text{ for every }t\ge0\}=1.
\end{equation}
Equations
\eqref{eq:final-global-capacity} and
\eqref{eq:final-noncollision} establish the assertions of
Theorem~\ref{thm:noncollision}.

\appendix
\section{Uniform GMC estimates}\label{sec:gmc}

\subsection{Scaled masses and Green energies}

Let $Q_0=(-1/2,1/2)^2$ and let $D_0$ be a fixed smooth
domain with $\overline{Q_0}\subset D_0$.

\begin{lemma}\label{lem:scaled-energy}
For every $s>0$,
\begin{equation}\label{eq:scaled-energy}
 \sup_{Q\in W}\E[M_Q^{-s}]<\infty,
 \qquad
 \sup_{Q\in W}\E[A_Q^s]<\infty.
\end{equation}
\end{lemma}

\begin{proof}
The rescaled field on $Q_0$ has covariance
\[
 K_Q(u,v)=G_D(c_Q+r_Qu,c_Q+r_Qv).
\]
The two-dimensional Green upper bound and \eqref{eq:whitney-geometry}
show that, for the reference domain $D_0$,
there is $C_{\mathrm{cov}}<\infty$ such that
\begin{equation}\label{eq:fixed-covariance-comparison}
 K_Q(u,v)\le G_{D_0}(u,v)+C_{\mathrm{cov}}
 \quad(u,v\in Q_0,\ Q\in W).
\end{equation}
The right side is the covariance of the restriction of a fixed
Dirichlet GFF to $Q_0$ plus an independent spatially constant normal
variable of variance $C_{\mathrm{cov}}$.  The rescaled measure
$\widetilde\mu_Q(B)=r_Q^{-2}\mu(c_Q+r_QB)$ has the Wick chaos law
associated with $K_Q$.  Let $M_0$ be the GMC measure of this reference GFF on $Q_0$.
Kahane's inequality applied to the convex
function $(x+\varepsilon)^{-s}$, followed by
$\varepsilon\downarrow0$, gives
\begin{equation}\label{eq:uniform-negative-moments}
 \E M_Q^{-s}\le
 \exp\!\left(\frac{\gamma^2C_{\mathrm{cov}}}{2}(s^2+s)\right)
 \E[M_0(Q_0)^{-s}],
\end{equation}
The right-hand side is finite for every $s>0$ by the subcritical GMC
negative-moment theorem \cite[Theorem 2.12]{RV}: a closed Euclidean
ball of positive radius contained in the interior of $Q_0$ has finite
negative chaos-mass moments, and its mass is bounded above by $M_0(Q_0)$.
For $1<b<q_*$, the same covariance comparison, applied to $x^b$,
gives
\begin{equation}\label{eq:uniform-positive-moments}
 \sup_{Q\in W}\E M_Q^b
 \le \exp\!\left(\frac{\gamma^2C_{\mathrm{cov}}}{2}(b^2-b)\right)
       \E[M_0(Q_0)^b]<\infty.
\end{equation}
Choose $b$ with $\zeta(b)>2$.  Applying this comparison to
$B(u,r)\cap Q_0$ gives, uniformly in $Q$ and $u\in Q_0$,
\begin{equation}\label{eq:scaled-ball-moment}
 \E[\widetilde\mu_Q(B(u,r))^b]
 \le C\E[M_0(B(u,r)\cap Q_0)^b]
 \le Cr^{\zeta(b)},\qquad 0<r\le1.
\end{equation}
Choose $0<\alpha_0<(\zeta(b)-2)/b$.  Let $J_j$ be a
deterministic grid of $O(2^{2j})$ balls of radius $C2^{-j}$ covering
the reference square at scale $2^{-j}$.  Define
\[
 C_Q=C\left(1+M_Q^b+
 \sum_{j\ge0}\sum_{B\in J_j}
 2^{j\alpha_0b}\widetilde\mu_Q(B)^b\right)^{1/b}.
\]
Grid covering and \eqref{eq:scaled-ball-moment} show
\begin{equation}\label{eq:scaled-frostman}
 \widetilde\mu_Q(B(u,r))\le C_Qr^{\alpha_0},\qquad
 \sup_Q\E C_Q^b
 \le C+C\sum_{j\ge0}2^{-j[\zeta(b)-2-\alpha_0b]}<\infty.
\end{equation}
Since $G_D(c_Q+r_Qu,c_Q+r_Qv)\le C+\log^+(1/|u-v|)$,
\begin{align}
 I_G(\mu|_Q)
 &\le Cm_Q^2+C m_Q^2\int_0^1
 \min\{1,(C_Q/M_Q)r^{\alpha_0}\}\,\frac{dr}{r}\notag\\
 &\le Cm_Q^2[1+\log C_Q+|\log M_Q|].\label{eq:scaled-energy-pointwise}
\end{align}
The moments of $\log C_Q$ of all finite orders follow from
\eqref{eq:scaled-frostman}.  Those of $|\log M_Q|$ follow from
\eqref{eq:uniform-negative-moments} and any uniform positive moment
of $M_Q$ below $q_*$.  The second assertion follows.
\end{proof}

\subsection{Joint moments for separated squares}

\begin{lemma}\label{lem:separated-moments}
For $1<q<q_*$ and separated $Q,R\in W$,
\begin{equation}\label{eq:separated-moments}
 \E[\mu(Q)^q\mu(R)^q]\le C_q(r_Qr_R)^{2q}.
\end{equation}
\end{lemma}

\begin{proof}
By \eqref{eq:fixed-covariance-comparison} on each square and
\eqref{eq:green-comparison} across the two separated squares, enlarge $C_{\mathrm{cov}}$ if necessary, independently of $Q,R$,
so that the covariance
\[
 G^+(x,z)=
 \begin{cases}
 G_{D_0}\bigl((x-c_Q)/r_Q,(z-c_Q)/r_Q\bigr)+C_{\mathrm{cov}},
       &x,z\in Q,\\
 G_{D_0}\bigl((x-c_R)/r_R,(z-c_R)/r_R\bigr)+C_{\mathrm{cov}},
       &x,z\in R,\\
 C_{\mathrm{cov}}, &x\in Q,\ z\in R\text{ or }x\in R,\ z\in Q
 \end{cases}
\]
satisfies $G_D(x,z)\le G^+(x,z)$ on $Q\cup R$.
It is the covariance of two independent rescaled reference GFFs
plus a common independent $N(0,C_{\mathrm{cov}})$ variable.
Let $X,X^+$ be independent fields with covariances $G_D,G^+$,
respectively, and write $X_t=\sqrt{1-t}X+\sqrt t X^+$.
After nonnegative mollification within the respective squares, write
$G_D^{(\delta)},G^{+(\delta)}$ for the covariances and
$M_{A,\delta}(t)$ for the Wick-chaos masses of $A\in\{Q,R\}$,
with associated random measure $M_{t,\delta}$.
Since $G^{+(\delta)}-G_D^{(\delta)}\ge0$, Gaussian
integration by parts, with
$F_\varepsilon(a,b)=(a+\varepsilon)^q(b+\varepsilon)^q$, gives
\begin{align*}
 &\frac{d}{dt}\E F_\varepsilon(M_{Q,\delta}(t),M_{R,\delta}(t))\\
 &\quad=\frac{\gamma^2}{2}\sum_{A,B\in\{Q,R\}}
 \iint_{A\times B}(G^{+(\delta)}-G_D^{(\delta)})(x,z)\\
 &\qquad\quad\cdot\E\bigl[\partial_{AB}F_\varepsilon
       (M_{Q,\delta}(t),M_{R,\delta}(t))
       M_{t,\delta}(dx)M_{t,\delta}(dz)\bigr]\ge0.
\end{align*}
Indeed, for $a,b\ge0$ and $q>1$,
\[
 \partial_{aa}F_\varepsilon
  =q(q-1)(a+\varepsilon)^{q-2}(b+\varepsilon)^q\ge0,
 \quad
 \partial_{ab}F_\varepsilon
  =q^2(a+\varepsilon)^{q-1}(b+\varepsilon)^{q-1}\ge0,
\]
and likewise $\partial_{bb}F_\varepsilon\ge0$.
Writing $M_A(1)$ for the limiting reference chaos mass, GMC
convergence and Fatou's lemma yield
\begin{align*}
 \E F_\varepsilon(\mu(Q),\mu(R))
 &\le\liminf_{\delta\downarrow0}\E F_\varepsilon
      (M_{Q,\delta}(0),M_{R,\delta}(0))\\
 &\le\lim_{\delta\downarrow0}\E F_\varepsilon
      (M_{Q,\delta}(1),M_{R,\delta}(1))
 =\E F_\varepsilon(M_Q(1),M_R(1)).
\end{align*}
The last limit uses the $q$-moment convergence of the reference
chaoses for $q<q_*$.  Let $\mu_Q^{(0)}$ and $\mu_R^{(0)}$ be the
independent rescaled reference chaoses without the common Gaussian
variable.  Their $q$-moments satisfy
\[
 \E[\mu_Q^{(0)}(Q)^q]\le Cr_Q^{2q},\qquad
 \E[\mu_R^{(0)}(R)^q]\le Cr_R^{2q}.
\]
Letting $\varepsilon\downarrow0$ and integrating the common Gaussian
variable yields
\begin{align*}
 \E[\mu(Q)^q\mu(R)^q]
 &\le \exp\{C_{\mathrm{cov}}\gamma^2(2q^2-q)\}
    \E[\mu_Q^{(0)}(Q)^q]\E[\mu_R^{(0)}(R)^q]\\
 &\le C_q(r_Qr_R)^{2q}.
\end{align*}
\end{proof}

\section{Uniform estimates along Gaussian circles}\label{sec:rotation-technical}

\begin{proof}[Proof of Lemma~\ref{lem:rotation-pair}]
Write $m_Q(\theta)=\mu_\theta(Q)$,
$m_R(\theta)=\mu_\theta(R)$, and
$F(\theta)=m_Q(\theta)m_R(\theta)$.  By
\eqref{eq:rotated-derivative},
\[
 F'(\theta)=\gamma\left[m_R(\theta)\int_Q\dot h_\theta\,d\mu_\theta
 +m_Q(\theta)\int_R\dot h_\theta\,d\mu_\theta\right].
\]
Conditional on $h_\theta$, this is a centered Gaussian variable.
Since the covariance $G_\tau$ is pointwise dominated by $G_D$,
\begin{equation}\label{eq:rotation-derivative-moment}
 \E[|F'(\theta)|^q\mid h_\theta]
 \le C_qF(\theta)^q
 \left[A_Q(\theta)^{1/2}+A_R(\theta)^{1/2}\right]^q.
\end{equation}
Here $A_Q(\theta)$ is the normalized Green energy
\eqref{eq:scaled-mass} of $\mu_\theta|_Q$.
Choose $q<q'<q_*$ and apply H\"older's inequality:
\begin{align*}
 \E[F^qA_Q^{q/2}]
 &\le(\E F^{q'})^{q/q'}
 \left(\E A_Q^{\frac{qq'}{2(q'-q)}}\right)^{1-q/q'}
 \le C(r_Qr_R)^{2q}.
\end{align*}
Lemma~\ref{lem:scaled-energy} controls the second factor and
Lemma~\ref{lem:separated-moments} the first.  Identical estimates hold
for $A_R$.  Also $\E|F|^q\le C(r_Qr_R)^{2q}$.
The circle Sobolev inequality
\[
 \sup_\theta|F(\theta)|^q
 \le C_q\int_0^{2\pi}(|F(\theta)|^q+|F'(\theta)|^q)\,d\theta
\]
and Fubini's theorem establish \eqref{eq:rotation-pair}.
\end{proof}

\begin{proof}[Proof of Lemma~\ref{lem:rotation-frostman}]
Set $b=a_0$.  Since $\zeta(b)>2+\alpha b$, choose $b'>b$ and
$\eta>0$ sufficiently small that
\begin{equation}\label{eq:beta-choice}
 b-\eta>1,\quad b'<q_*,\qquad
 \beta:=\min\left\{\zeta(b),\frac b{b'}\zeta(b'),
 \zeta(b-\eta)\right\}>2+\alpha b.
\end{equation}
Fix a ball $B=B(x,r)\cap D$ and set $m=\mu(B)$.
From \eqref{eq:frostman} and the Green kernel logarithmic bound,
\begin{equation}\label{eq:ball-green-mass}
 I_G(\mu|_B)\le Cm^2
 \left(1+\log C_\alpha+\log_+(1/m)\right).
\end{equation}
By H\"older's inequality and the moments of $\log C_\alpha$,
\begin{align}
 \E[m^b(1+\log C_\alpha)^{b/2}]
 &\le(\E m^{b'})^{b/b'}
    \left(\E(1+\log C_\alpha)^{\frac{bb'}{2(b'-b)}}
    \right)^{1-b/b'}\notag\\
 &\le Cr^{\frac b{b'}\zeta(b')},\label{eq:ball-green-holder}\\
 \E[m^b\log_+(1/m)^{b/2}]
 &\le C\E[m^{b-\eta}]
 \le Cr^{\zeta(b-\eta)}.\label{eq:ball-green-log}
\end{align}
Together with \eqref{eq:ball-moment}, these inequalities imply
\begin{equation}\label{eq:ball-energy-moment}
 \E[I_G(\mu|_B)^{b/2}]\le Cr^\beta.
\end{equation}
For $m_\theta(B)=\mu_\theta(B)$, independence of
$\dot h_\theta$ and $h_\theta$ gives
\begin{align*}
 \E|m_\theta'(B)|^b
 &\le C_b\E[I_{G_\tau}(\mu_\theta|_B)^{b/2}]\\
 &\le C_b\E[I_G(\mu_\theta|_B)^{b/2}]
 \le Cr^\beta.
\end{align*}
Using $\E m_\theta(B)^b\le Cr^{\zeta(b)}$, the circle Sobolev
inequality yields
\begin{equation}\label{eq:ball-rotation-sup}
 \E\left[\sup_\theta\mu_\theta(B)^b\right]
 \le C_b\int_0^{2\pi}
 \E\bigl[m_\theta(B)^b+|m_\theta'(B)|^b\bigr]d\theta
 \le Cr^\beta.
\end{equation}
For $r_j=2^{-j}$, choose a deterministic covering by $O(r_j^{-2})$
balls of radius comparable to $r_j$.  Define
\[
 C_\tau=C\left(1+\sum_{j\ge0}\sum_{B\text{ in layer }j}
 r_j^{-\alpha b}\sup_\theta\mu_\theta(B)^b\right)^{1/b}.
\]
The covering ensures the asserted Frostman inequality; and
\[
 \sup_\tau\E C_\tau^b
 \le C+C\sum_jr_j^{\beta-2-\alpha b}<\infty.
\]
The coarse-scale term ensures $C_\tau\ge1+\sup_\theta\mu_\theta(D)$
by modifying the harmless multiplicative constant.
\end{proof}

\end{document}